\documentclass[journal]{IEEEtran}

\ifCLASSINFOpdf
  \usepackage[pdftex]{graphicx}
\else
\fi

\usepackage{cite}
\usepackage{amsmath,amssymb,amsfonts}
\usepackage{algorithm}
\usepackage{algpseudocode}
\usepackage{graphicx}
\usepackage{tabularx}
\usepackage{textcomp}
\usepackage{float}
\usepackage{mathtools,nccmath}
\usepackage{multirow}
\usepackage{booktabs}
\usepackage{longtable}

\ifCLASSOPTIONcompsoc
    \usepackage[caption=false, font=normalsize, labelfont=sf, textfont=sf]{subfig}
\else
\usepackage[caption=false, font=footnotesize]{subfig}
\fi

\algdef{SE}[SUBALG]{Indent}{EndIndent}{}{\algorithmicend\ }%
\algtext*{Indent}
\algtext*{EndIndent}

\def\BibTeX{{\rm B\kern-.05em{\sc i\kern-.025em b}\kern-.08em
    T\kern-.1667em\lower.7ex\hbox{E}\kern-.125emX}}

\begin{document}

\title{Efficient Real-time Rail Traffic Optimization: Decomposition of Rerouting, Reordering, and Rescheduling Problem}
\author{László Lindenmaier, István Ferenc Lövétei, Szilárd Aradi
\thanks{Department of Control for Transportation and Vehicle Systems, Budapest University of Technology and Economics, 1111 Budapest, Hungary e-mail: \{lindenmaier.laszlo, lovetei.istvan, aradi.szilard\}@kjk.bme.hu}
\thanks{This work was supported in part by the Ministry of Innovation and Technology NRDI Office within the
framework of the Autonomous Systems National Laboratory Program; in part by the Ministry of Innovation and Technology of Hungary from the National Research, Development and Innovation Fund, under the TKP2021 Funding Scheme, under Project BME-NVA-02.

This paper was supported by the János Bolyai Research Scholarship of the Hungarian Academy of Sciences.}}

\maketitle

\begin{abstract}
The railway timetables are designed in an optimal manner to maximize the capacity usage of the infrastructure concerning different objectives besides avoiding conflicts. The real-time railway traffic management problem occurs when the pre-planned timetable cannot be fulfilled due to various disturbances; therefore, the trains must be rerouted, reordered, and rescheduled. Optimizing the real-time railway traffic management aims to resolve the conflicts minimizing the delay propagation or even the energy consumption. In this paper, the existing mixed-integer linear programming optimization models are extended considering a safety-relevant issue of railway traffic management, the overlaps. However, solving the resulting model can be time-consuming in complex control areas and traffic situations involving many trains. Therefore, we propose different runtime efficient multi-stage heuristic models by decomposing the original problem. The impact of the model decomposition is investigated mathematically and experimentally in different rail networks and various simulated traffic scenarios concerning the objective value and the computational demand of the optimization. Besides providing a more realistic solution for the traffic management problem, the proposed multi-stage models significantly decrease the optimization runtime.
\end{abstract}


\section{Introduction}
\label{sec:introduction}
\IEEEPARstart{T}{he} low specific energy consumption of rail transportation compared to the other sectors is great potential to reduce CO\textsubscript{2} emission \cite{GARCIAOLIVARES2018266, DALLACHIARA2017227, WANG2021270, su6020702}. The climate protection initiatives and the increasing energy prices will enable the railway to play a greater role in the transport of the future, increasing the demand for passenger and freight rail services \cite{chaturvedi2015long, islam2015assessing, WANG2014641}. Moreover, increasing infrastructure utilization, competitiveness and decreasing energy consumption and carbon emission is one of the most crucial research and development areas in railway traffic \cite{GONZALEZGIL2014509, LI2020121316, MARTINEZFERNANDEZ2019153}. The timetable optimization is one of the key components of increasing energy efficiency \cite{7151801, en13051115}. However, in many cases, the optimally scheduled trains suffer primary delays due to different traffic disturbances, and the pre-planned timetable would result in conflicting trains \cite{Nagy_Csiszar_2015}. Therefore, the trains must be rescheduled or rerouted and reordered, assigning secondary delays to the trains to solve the real-time railway traffic management problem (rtRTMP) \cite{corman2010real}. Nowadays, the conflicts are resolved by human dispatchers according to several generic and regional principles, such as the first-in-first-out (FIFO) policy. However, making an optimal decision within the required short response time is very challenging for a human dispatcher. Hence, automatic railway traffic management and control is an emerging research topic aiming to reduce energy consumption and optimize the capacity usage of rail transportation to satisfy the increasing demands \cite{PELLEGRINI201613, botte2018dispatching}.

\subsection{Related Work}
Many works tackle the real-time railway traffic management problem with different optimization approaches and techniques to minimize delays or even energy consumption \cite{6920082, CACCHIANI201415}. The models can be divided into microscopic and macroscopic groups according to the levels of the infrastructure representation \cite{CACCHIANI201415, 10.1007/978-3-319-24264-4_45}. The microscopic models include detailed information about the rail network in graphs, as in \cite{RODRIGUEZ2007231} and \cite{d2008reordering}. The macroscopic models neglect some lower-level resources such as signals, reducing the complexity of the infrastructure model \cite{dollevoet2015delay, kanai2011optimal}. 

The optimization algorithms for the rtRTMP range widely. Mixed-integer linear programming (MILP) is one of the commonly used frameworks, as in \cite{acuna2011mip, tornquist2015configuration}. Pellegrini et al. proposed a MILP model in \cite{pellegrini_2012, PELLEGRINI201458}, considering the infrastructure at a microscopic, track circuit level. A macroscopic-microscopic approach decomposes the problem into line traffic and station control problems in \cite{manino2011real, lamorgese2015exact, lamorgese2013track}. However, since the problem is NP-hard in general, some custom heuristics are introduced in \cite{boccia2013dispatching, pellegrini2015recife} to decrease the computational time of the algorithms. Furthermore, Pellegrini et al. proposed an efficient implementation in \cite{pellegrini2019efficient} by extending their previous work presented in \cite{pellegrini2015recife} with new valid constraints.

Besides MILP formulation, many other optimization algorithms exist to solve the rtRTMP. A distributed model predictive control (DMPC) is proposed in \cite{kersbergen2016distributed} with similar mixed-integer linear constraints, as in \cite{pellegrini_2012}. Ant colony optimization (ACO) is applied in \cite{SAMA201689} to find the best subset of routing alternatives that serve as an input for the original problem formulated as mixed-integer linear programming. While the differences between ACO and MILP regarding the preliminary route selection are pointed out in \cite{PASCARIU2021167}. The lower level control of automatic train operation (ATO) is connected with the traffic management optimization in \cite{rao2013holistic, rao2016new}. While the adaptive train control and speed profile optimization is integrated within an MPC framework in \cite{caimi2012model}. Törnquist proposed an effective greedy approach in \cite{TORNQUISTKRASEMANN201262} to satisfy the requirement of the available decision time. Other soft computing techniques with lower computational requirements, such as Monte Carlo tree search and reinforcement learning (RL), are used as alternatives to conventional approaches \cite{lovetei2021mcts, lovetei2022environment}.

The different models can also be distinguished according to the objective function. In most of the listed researches (e.g., \cite{pellegrini2015recife, pellegrini2019efficient, SAMA201689}), the optimization goal is to to reduce the delay propagation concerning all trains. However, some papers focus on energy consumption, as in \cite{naldini:hal-03467003, naldini2021ant, MONTRONE2018524}, or a trade-off is made between time and energy efficiency \cite{rao2013railway, d2010running, 7057643}. Some works aim to increase the comfort of the passengers by minimizing the platform changes and passenger train delays, as in \cite{dollevoet2012delay, corman2015railway}, compared to the ones that do not distinguish between freight and passenger trains \cite{CAIMI20122578, dundar2013train}. Another aspect is to make a trade-off between the delay of all trains reduction and passenger satisfaction by cancelling some transfer connections \cite{corman2012bi, articleGinkel}.

\subsection{Contributions of the Paper}
Mixed-integer linear programming (MILP) allows considering railway traffic regulations by constraining the optimization space with straightforward expressions. Moreover, it provides a global optimum for the problem if the response time is not crucial. Therefore, MILP is a commonly used algorithm to solve the rtRTMP, as shown in the previous section. However, since the problem is NP-hard in general, in the case of large networks and complex scenarios, the runtime of the optimization can exceed the available time. Moreover, we extend the existing models in \cite{pellegrini_2012, pellegrini2015recife} to consider the safety concept of overlaps. Since the model extension results in high computational cost, a runtime-efficient implementation is essential for optimizing the rtRTMP in a real traffic environment. Some works aim to decrease the computational complexity of MILP by introducing custom heuristics and additional constraints, as in \cite{pellegrini2015recife, pellegrini2019efficient}. Preliminary route selection is an alternative approach to accelerating optimization by reducing the number of possible solutions, as in \cite{SAMA201689}. This paper proposes a multi-stage framework consisting of different MILP models, minimizing the delay propagation by decomposing the rtRTMP into rerouting, reordering, and rescheduling sub-problems. A low-complexity model is implemented to formulate the objective function that focuses on the relation between possible train pairs. This simplified approach provides a fast, sub-optimal solution for the rtRTMP by utilizing the results as constraints. The result of the simplified model can also serve as an initial condition for the second stage of optimization, reducing the response time significantly when finding the first feasible solution would last longer. Furthermore, we analyzed the optimality condition to obtain the circumstance in which this solution results in a sub-optimal solution. As a result of optimality condition analysis, we developed a new model that, besides decreasing the computational cost, it guarantees the global optimum for the original problem. Similar to the simplified approach, the runtime-efficient global optimum model supports the second stage of the optimization solving the extended problem. The optimization model is introduced in Section \ref{section:Model Formulation} with a brief overview of the existing solutions. The proposed sub-optimal solution is detailed in Section \ref{section:suboptimal model}. The optimality condition of the first approach is analyzed in Section \ref{section:optimality condition}, and the resulting global optimum model is explained in Section \ref{section:global optimum model}. Finally, the results evaluated in a simulation environment are shown in Section \ref{section:results}.

\section{Model Formulation} \label{section:Model Formulation}

The proposed optimization algorithm is formulated as a mixed-integer linear programming problem in line with previous works in \cite{pellegrini_2012, pellegrini2015recife} that minimizes the delay propagation within a control area. We extended the existing solutions to consider a safety-relevant traffic management regulation, the overlaps. Although the model extension contributes to more feasible solutions in practice, it increases the computational cost handled by the model decomposition. The rtRTMP to be solved is defined by the trains assigned in space and time to the control area. The multi-layer infrastructure model is constructed from a detailed microscopic graph representation of the rail network, as in \cite{lindenmaier2021infrastructure}. The train schedule is given concerning the fundamental elements of the railway infrastructure, the track circuits represented by the corresponding edges of the graph model, on which the presence of a train is automatically detected. The group of track circuits delimited by signals form block sections whose consecutive order compose different routes between the boundaries of the control area. The following nomenclature in accordance with \cite{pellegrini_2012, pellegrini2015recife} is used for the inputs of the optimization:

\begin{table}[htbp]
\setlength\extrarowheight{2.5pt}
\begin{center}
\begin{tabularx}{\linewidth}{lX}
$TC$, $R$ & set of track circuits and routes of control area \\
$TC_r$ & track circuits constructing route $r$ \\
$tc_0$, $tc_{\infty}$ & dummy track circuits denoting the entry and exit locations of the control area \\
$p_{r,tc}$, $s_{r,tc}$ & preceding and subsequent track circuits of $tc$ along route $r$\\
$P_{r,tc}$ & set of track circuits preceding $tc$ along route $r$\\
$sw_{r,tc}^o$, $sw_{r,tc}^c$ & binary indicators with a value of $1$ if track circuit $tc$ involves an opening or closing switch along $r$, respectively \\
$bs_{r,tc}$ & block section of track circuit $tc$ on route $r$\\
$ref_{r,tc}$ & reference track circuit to reserve $tc$ on route $r$ (specifically the first track circuit of $bs_{r,tc}$ in a two-state signaling system) \\
$T$ & set of trains assigned to the control area \\
$TC^t \subseteq TC$, $R^t \subseteq R$ & set of track circuits and routes available for train $t$ \\
$w^t$ & priority weight of train $t \in T$ \\
$init^t$, $exit^t$ & the time when train $t$ is scheduled to enter and leave the control area \\
$tc_{in}^t$, $tc_{ex}^t$ & the track circuits at where train $t$ enters and leaves the control area \\
$ds_{r,tc}^t$, $cl_{r,tc}^t$ & total duration of stay and clearing time of train $t$ on track circuit $tc$ along route $r$ \\
$for$, $rel$, $over$ & formation, release, and overlap time \\
$M$ & large constant
\end{tabularx}
\end{center}
\label{table: optimization inputs}
\vspace{-2mm}
\end{table}

\noindent The optimization state space $X$ is formed by continuous variables and binary indicators with the following notations:

\begin{table}[htbp]
\setlength\extrarowheight{2.5pt}
\begin{center}
\begin{tabularx}{\linewidth}{lX}
$e_{r, tc}^t$ & time when train $t$ enters track circuit $tc$ along $r$ \\
$d_{r, tc}^t$ & delay assigned to train $t$ at track circuit $tc$ on route $r$ \\
$sR_{tc}^t$, $eR_{tc}^t$ & time when train $t$ starts reserving and releases track circuit $tc$ \\
$sO_{tc}^t$, $eO_{tc}^t$ & time when train $t$ starts and finishes reserving track circuit $tc$ due to overlap\\
$x_r^t$ & route indicator with a value of $1$ if train $t$ travels along route $r$, $0$ otherwise \\
$z_{tc}^t$ & stop indicator with a value of $1$ if $over$ allocated for the safe stopping expires before train $t$ can utilize track circuit $tc$, $0$ otherwise \\
$c_{t', tc}^t$ & conflict indicator with a value of $1$ if train $t$ and $t'$ would collide at track circuit $tc$ without rescheduling \\
$y_{t', tc}^t$ & precedence indicator with a value of $1$ if train $t$ utilizes track circuit $tc$ before $t'$ \\
$yO_{t', tc}^t$ & overlap precedence indicator with a value of $1$ if train $t$ reserves track circuit $tc$ due to overlap before $t'$ \\
$D^t$ & total secondary delay of train $t$ within the control area 
\end{tabularx}
\end{center}
\label{table: train state variables}
\end{table}
\vspace{-5mm}
\begin{algorithm}[H]
    \caption{Scheduling constraints}
    \label{alg:scheduling constraints}
    \begin{algorithmic}[1] 
        \For{$\forall t \in T$}
            \State $e_{tc_{\infty}}^t = 0$
            \For{$\forall r \in R^t$}
                \For{$\forall tc \in TC_r$}
                    \Statex Time constraints: 
                    \State $e_{r, tc}^t \leq M x_r^t$
                    \State $e_{r, tc}^t \geq e_{r, p_{r,tc}}^t + ds_{r,tc}^t \, x_r^t$
                    \If {$p_{r,tc} = tc_0$}
                        \State $e_{r, tc}^t = init^t \, x_r^t$
                    \EndIf
                    
                     \If {$s_{r,tc} = tc_{\infty}$}
                        \State $e_{tc_{\infty}}^t = e_{tc_{\infty}}^t + e_{r, tc}^t + ds_{r,tc}^t \, x_r^t$
                    \EndIf
                    
                    \Statex Delay constraints:
                    \If {$bs_{r,tc} \neq bs_{r,s_{r,tc}}$}
                        \State $d_{r, tc}^t = e_{r, s_{r,tc}}^t - e_{r, tc}^t - ds_{r,tc}^t$
                    \Else
                        \State $d_{r, tc}^t = 0$
                    \EndIf
                \EndFor
            \EndFor
            \State $D^t \geq e_{tc_{\infty}}^t - exit^t$
        \EndFor
    \end{algorithmic}
\end{algorithm}

\begin{algorithm}
    \caption{Capacity constraints}
    \label{alg:capacity constraints}
    \begin{algorithmic}[1] 
        \For{$\forall t \in T$}
            \State $\sum_{\forall r \in R^t} x_r^t = 1$ \label{constr: route}
            \For{$\forall tc \in TC^t$}
                \State $sR_{tc}^t = \sum \limits_{\substack{\forall r \in R^t: \\ tc \in TC_r}} e_{r, ref_{r,tc}}^t - for \: x_r^t$
                \vspace{1.5mm}
                \State $eR_{tc}^t = \sum \limits_{\substack{\forall r \in R^t: \\ tc \in TC_r}} e_{r, ref_{r,tc}}^t + (cl_{r,tc}^t + rel) \: x_r^t$
                \vspace{1.5mm}
                \For{$\forall t' \in T: t' \neq t \land tc \in TC^{t'}$}
                    \State $y_{t', tc}^t + y_{t, tc}^{t'} = 1$ 
                    \label{constr: precedence}
                    \State $sR_{tc}^t \geq eR_{tc}^{t'} - M \, y_{t', tc}^t$ \label{constr: conflict1}
                    \State $sR_{tc}^{t'} \geq eR_{tc}^{t} - M(1 - y_{t', tc}^t)$ \label{constr: conflict2}
                \EndFor
            \EndFor
        \EndFor
    \end{algorithmic}
\end{algorithm}

According to \cite{pellegrini_2012}, the model constraints can be distinguished into the time concerning, delay managing, rolling stock configuration related, and capacity constraints considering the railway traffic regulations. In our model, the scheduling constraints involving the time and delay constraints are detailed in Alg. \ref{alg:scheduling constraints}. The time constraints ensure the temporal coherence of train schedules and allow distinguishing the same track circuits belonging to different routes. The delay constraints manage the local delays assigned to the trains at track circuits considering the block sections determined by the signalling system. The $D^t$ total delay as the basis of the objective function is bounded by the difference between the actual and scheduled exit time of the train.

The capacity constraints in line \ref{constr: route} of Alg. \ref{alg:capacity constraints} impose that every train utilizes exactly one route. According to safety regulations, a train $t$ reserves all $tc$ track circuit of $bs_r$ block section with $for$ time in advance before it enters the first track circuit $ref_{r,tc}$ of $bs_r$ along route $r$. Moreover, the reservation of a track circuit $tc$ by train $t$ is released only if $tc$ is left by the complete vehicle assembly according to $cl_{r,tc}^t$ and the $rel$ release time expires. Finally, the constraints from line \ref{constr: precedence} to \ref{constr: conflict2} prevent track circuits from being reserved by multiple trains at a given time. Since the focus of this paper is the runtime efficiency, the trains in connection and the constraints due to the change of rolling stock configuration with less impact on the complexity are not considered.

We extend the existing model in Alg.~\ref{alg:overlap constraints}, considering that the trains may not be able to stop before signals due to the driver's reaction time or increased braking distance. In railway traffic control, a safety distance, referred to as overlap, is provided beyond stop signals for the trains forced to stop. According to the constraints in line \ref{constr: start overlap}, besides the used block section, the first track circuit following the signal that terminates the section is also reserved for the train. The track circuit is reserved due to overlap until the overlap time $over$ within what the train is assumed to be able to stop safely expires. However, if the train does not stop according to $z_{tc}^t$, the reservation due to overlap finished at time $eO_{tc}^t$ is bounded by $sR_{tc}^t$ in line \ref{constr: end overlap2}. The conflicting reservation of a track circuit is prevented by the constraints in lines \ref{constr: conflict} - \ref{constr: overlap conflicts}.
\vspace{1mm}
\begin{algorithm}
    \caption{Overlap constraints}
    \label{alg:overlap constraints}
    \begin{algorithmic}[1] 
        \For{$\forall t \in T$}
            \For{$\forall tc \in TC^t$}
            \State $eOver_{tc}^t = \sum_{\substack{\forall r \in R^t: \\ tc \in TC_r}} e_{r, p_{r,tc}}^t + over \: x_r^t$
            \State $z_{tc}^t \geq \dfrac{1}{M} \left(sR_{tc}^t - eOver_{tc}^t \right)$
                \If {$bs_{r,tc} \neq bs_{r,p_{r,tc}}$}
                    \State $sO_{tc}^t = \sum_{\substack{\forall r \in R^t: \\ tc \in TC_r}} e_{r, ref_{r,p_{r,tc}}}^t - for \: x_r^t$ \label{constr: start overlap}
                    \vspace{1mm}
                    \State $eOver_{tc}^t - M(1 - z_{tc}^t) \leq eO_{tc}^t \leq eOver_{tc}^t$ \label{constr: end overlap1}
                    \State $sR_{tc}^t - M \, z_{tc}^t \leq eO_{tc}^t \leq sR_{tc}^t$
                    \label{constr: end overlap2}
                \EndIf
                
                \For{$\forall t' \in T: t' \neq t \land tc \in TC^{t'}$}
                    \If {$bs_{r,tc} \neq bs_{r,p_{r,tc}}$}
                        \State $c_{t',tc}^t = c_{t,tc}^{t'}$ \label{constr: conflict}
                        \State $eO_{tc}^{t'} - M(2 - c_{t',tc}^t - y_{t',tc}^t) \leq sR_{tc}^t$
                        \State $eO_{tc}^{t} - M((1 - c_{t',tc}^t) + y_{t',tc}^t) \leq sR_{tc}^{t'}$
                        \State $eR_{tc}^{t} - M(c_{t',tc}^t + (1 - y_{t',tc}^t)) \leq sO_{tc}^{t'}$
                        \State $eR_{tc}^{t'} - M(c_{t',tc}^t + y_{t',tc}^t) \leq sO_{tc}^{t}$
                        \vspace{1mm}
                        \State $yO_{t', tc}^t + yO_{t, tc}^{t'} = 1$ 
                        \State $sO_{tc}^t \geq eO_{tc}^{t'} - M \, yO_{t', tc}^t$
                        \State $sO_{tc}^{t'} \geq eO_{tc}^{t} - M(1 - yO_{t', tc}^t)$ \label{constr: overlap conflicts}
                    \EndIf
                \EndFor
            \EndFor
        \EndFor
    \end{algorithmic}
\end{algorithm}

\section{Model Decomposition} \label{section:methodology}

\begin{figure}
    \centering
    \includegraphics[width=0.7\columnwidth]{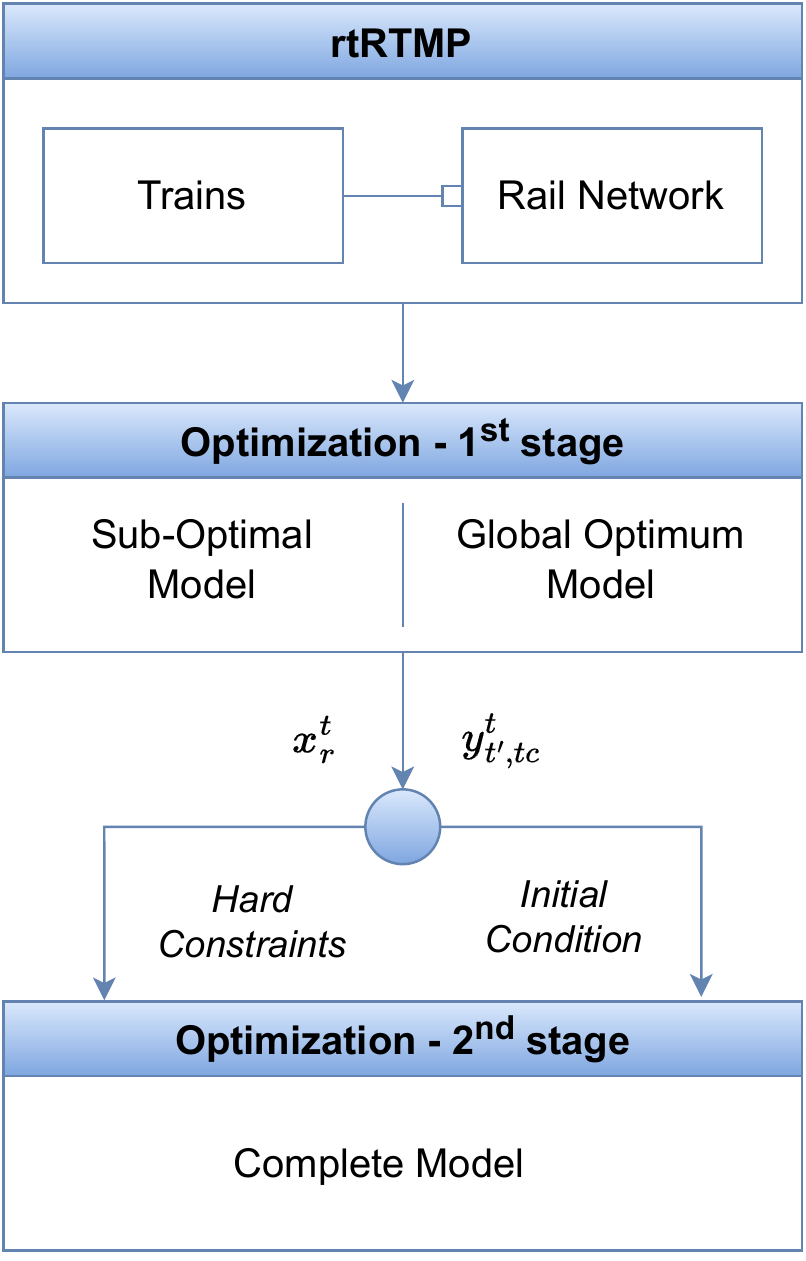}
    \caption{The flowchart of the proposed optimization model}
    \label{fig:flowchart}
\end{figure}

The complexity of the extended model consisting of many new integer variables and corresponding constraints causes a significant increase in runtime. Different heuristics and efficient implementations are accelerating the convergence of optimization as in \cite{boccia2013dispatching, pellegrini2015recife, pellegrini2019efficient}. The model decomposition allows dividing the original model into sub-problems with lower complexity as in \cite{SAMA201689}. The proposed decomposition is formulated based on the MILP framework described in Section \ref{section:Model Formulation} provides easy interoperability between the sub-models and the original problem. Therefore, the solutions of the sub-problems can support the complete MILP model to decrease the computational cost, forming a multi-stage optimization model. The results of the sub-models can form hard constraints in the original problem or serve as an initial condition for the second stage of optimization. Since the models in the first stage of the optimization neglect some constraints of the rtRTMP, only the $x_r^t$ route selection and the $y_{t',tc}^t$ precedence indicators are used as the solutions to the rerouting and reordering sub-problems. The flowchart of the proposed optimization model is illustrated in Fig.~\ref{fig:flowchart}.

\subsection{Sub-Optimal Model} \label{section:suboptimal model}
\begin{figure}
    \centering
    \subfloat[Case of trains traveling in opposite directions]{\includegraphics[width=0.99\linewidth]{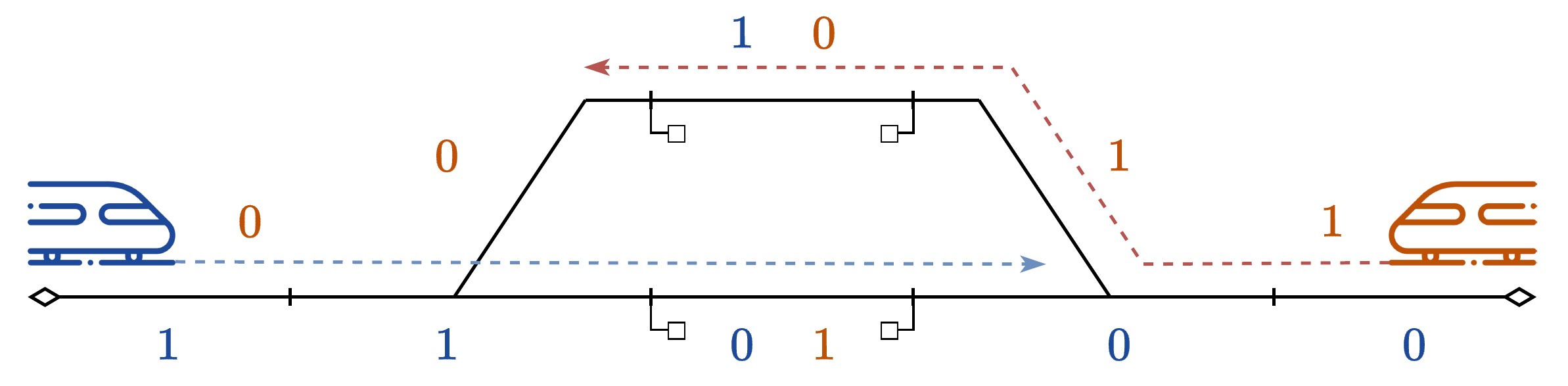} \label{fig:precedence oncoming}}
    
    \hfil
    
    \subfloat[Case of overtaking train] {\includegraphics[width=0.99\linewidth]{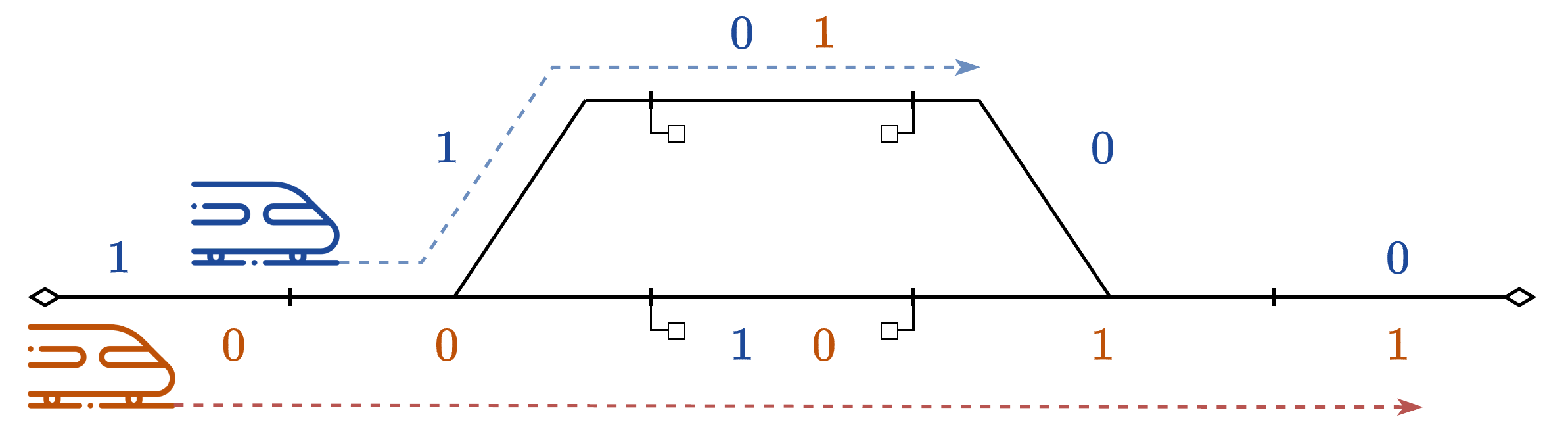} \label{fig:precedence overtaking}}
    \caption{The process of changing the precedence indicators}
    \label{fig:precedence change}
    \vspace{-3mm}
\end{figure}
The basic idea behind the sub-optimal model is to obtain a feasible solution for the rtRTMP close to the global minimum of train delays with a fast response time. Therefore, the original model is simplified by omitting the overlap constraints and reformulating the scheduling and capacity constraints. Delaying a train may affect all the other trains, thereby upsetting the entire schedule. Hence, the rtRTMP cannot be managed locally, unlocking the conflicts between the train pairs. However, the total secondary delay is related to the $d_{t',tc}^t$ delays interpreted between $\forall \{t,t'\} \in T$ pairs at track circuit $tc$ defined as:
\begin{equation}
    d_{t',tc}^t = 
    \begin{cases}
        \Hat{e}R_{tc}^{t} - \Hat{s}R_{tc}^{t'}, & 
        \text{if \ } y_{t',tc}^t > 0 \\
        \Hat{e}R_{tc}^{t'} - \Hat{s}R_{tc}^{t}, & \text{otherwise}
    \end{cases}
    \: , \label{eq:train pair delay}
\end{equation}
where $\Hat{s}R_{tc}^{t}$ and $\Hat{e}R_{tc}^{t}$ denote the times when train $t$ would start to reserve and release track circuit $tc$ without delays. According to \eqref{eq:train pair delay}, the delay resolving the conflict of train $t$ and $t'$ at $tc$ can be computed as the difference of the time when the reservation of the track circuit is finished by the train using it first, and when the other one would start the reservation according to the original timetable. Hence, the schedule variables are not updated in the sub-optimal model, neglecting the capacity constraints in lines \ref{constr: conflict1} - \ref{constr: conflict2} of Alg.~\ref{alg:capacity constraints}. Moreover, the delay constraints are reformulated defining the $d_{t',tc}^{t}$ delay assigned to train $t$ at track circuit $tc$ due to $t'$ as:
\begin{align}
    & d_{t',tc}^{t} \geq eR_{tc}^{t'} - sR_{tc}^{t} - \!\!\! \sum_{\substack{\forall r \in R^t: \\ tc \in TC_r}} \, \sum_{\substack{\forall p_{tc} \in \\ P_{r,tc}}} d_{t',p_{tc}}^t \!\! - M \, y_{t',tc}^t \, , \label{eq:submodel delay1} \\
    & d_{t,tc}^{t'} \geq eR_{tc}^{t} - sR_{tc}^{t'} - \!\!\! \sum_{\substack{ \forall r \in R^{t'}: \\ tc \in TC_r}} \, \sum_{\substack{\forall p_{tc} \in \\ P_{r,tc}}} d_{t,p_{tc}}^{t'} \!\! - M \, y_{t,tc}^{t'} \, , \label{eq:submodel delay2}
\end{align}
where the nested sums consider the already served waiting times, preventing multiple delays. For sake of brevity, the non-negativity boundaries of the delays are not given explicitly. Then, the objective function of the sub-optimal model considering the possible train pairs is formed by the reshaped $d_{t',tc}^{t}$ delay triplets as:
\begin{equation}
    f_s \left(d_{t',tc}^{t}, \delta^t \, \rvert \, X \right) = \!\! \sum_{\forall t \in T} w^t \! \left( \sum_{\substack{\forall t' \in T: \\ t' \neq t}} \!\!\! \sum_{\:\:\:\: \forall tc \in TC^t} \!\!\!\!\! d_{t',tc}^{t} \!+ \delta^t \!\! \right),
\end{equation}
where the $\delta^t$ additional delay resulting from alternative route selection is computed as:
\begin{equation}
    \delta^t = \sum_{\forall r \in R^t} \left( e_{r, tc}^t + ds_{r,tc}^t \, x_r^t \right) - exit^t, \quad s_{r,tc} = tc_{\infty}
\end{equation}

The simplified model must provide a feasible solution so that the result can be used in the second stage of the optimization. Therefore, additional constraints related to the precedence indicators are imposed in Alg.~\ref{alg:conflict constraints}. The process of changing the precedence indicators between two trains is illustrated in Fig.~\ref{fig:precedence change}, showing that the order of the trains can change at track circuits $tc$ that include a switch. However, the $y_{t',tc}^t$ precedence indicator may also change when the route of the trains splits at the preceding track circuit $p_{r,tc}$. The constraint in line \ref{constr: prec3} of Alg.~\ref{alg:conflict constraints} prevents changing $y_{t',tc}^t$ when $tc$ or $p_{r,tc}$ does not contain a switch. Moreover, the precedence indicator $y_{t',tc}^t$ cannot be incremented to $1$ only if $tc$ is not used by the $t$ or the overtaken train $t'$ does not use the track circuit preceding $tc$ along the route travelled by $t$ as in Fig.~\ref{fig:precedence overtaking}. Reversing this constraint,  $y_{t',tc}^t$ may change to $0$ from $1$ when $tc$ is not used by the other train $t'$ or the preceding track circuit of $t$ is not used by $t'$. The expressions from lines \ref{constr: prec4} to \ref{constr: prec7} impose these constraints, where $sw_{r,tc}^o$ means that the switch branches out the route, while $tc$ with $sw_{r,tc}^c = 1$ closes it, having less subsequent track circuits than preceding ones. According to the constraints in lines \ref{constr: prec8}~-~\ref{constr: prec9}, train $t''$ cannot use track circuit $tc$ before $t$ if it enters $tc$ before $t'$ and $t'$ before $t''$. The precedence $y_{t',tc}^t$ between train $t$ and $t'$ at track circuit $tc$ is set to $0$ by inequalities in lines \ref{constr: prec10}~-~\ref{constr: prec11} if $t$ travels along $tc$ while $t'$ does not.

However, besides the constraints in Alg.~\ref{alg:conflict constraints}, additional inequalities imposed $\forall \{t, t'\} \in T$ train pairs at their common $TC^t \cap TC^{t'}$ track circuits are needed to ensure that the trains enter the control area when they are scheduled. Since the sub-optimal model does not consider the signalling system, the delay $d_{t',tc}^{t}$ resolving the conflict between trains $t$ and $t'$ at track circuit $tc$ would be assigned to $t$ at the block section preceding $bs_{r,tc}$. Therefore, delaying train $t$ outside the control area is prevented by setting $d_{t',tc}^{t}$ to zero if $tc$ belongs to the first block section of $t$'s route as:
\begin{equation}
     d_{t',tc}^{t} \leq M \!\!\!\! \sum_{\substack{r \in R^{t}: \\ tc \in TC_r}} \!\!\!\!\! x_{r}^t \, ebs_{r,tc}^{t} \, , \label{eq:submodel delay3}
\end{equation}
where $ebs_{r,tc}^t \in \{0,1\}$ indicates if the block section $bs_{r,tc}$ of track circuit $tc$ is the first one along route $r$ of train $t$ as:
\begin{equation}
     ebs_{r,tc}^{t} = \begin{cases}
         1 & \text{if \ } bs_{r,tc} = bs_{r,tc_{in}^t} \\
         0 & \text{otherwise}
     \end{cases}
\end{equation}
However, the delay $d_{t',tc}^{t}$ formulated by $t$ and $t'$ neglects the affect of other trains. Supposing train $t'$ leaves the control area before $t$ enters it at the same track circuit $tc_{in}^t = tc_{ex}^{t'}$ according to the original timetable, the delays assigned to $t'$ due to another $t'' \in T$ may change the precedence of $t$ and $t'$. Therefore, the delay accumulated by $t'$ up to the block section where $t$ enter the control area is tackled by the constraint imposed in \eqref{eq:submodel delay4}. Moreover, delaying a train $t'$ traveling in the same direction as $t$ at the track circuit $tc_{in}^t$ forces $t$ to enter control area later then it is scheduled according to the timetable. Therefore, the $d_{t'',s_{tc}}^{t'}$ delays assigned to $t'$ at track circuit $tc$ due to the conflict with with another $t''$ at track circuits $s_{tc}$ is prevented by the expression in \eqref{eq:submodel delay5} if $tc$ belongs to the first block section of $t$ and $t'$ uses $tc$ before $t$. For sake of brevity, the counter-pair of the constraints in \eqref{eq:submodel delay3}, \eqref{eq:submodel delay4}, \eqref{eq:submodel delay5} related to the $d_{t,tc}^{t'}$ delay resolving the conflict of $t$ and $t'$, and the $d_{t'',p_{tc}}^{t'}$, $d_{t'',s_{tc}}^{t'}$ propagated delays are neglected.
\begin{multline}
        sR_{tc}^{t} \geq eR_{tc}^{t'} + \sum_{\substack{\forall t'' \in T: \\ t'' \neq t}} \, \sum_{\substack{\forall r \in R^{t'}: \\ tc \in TC_r}} \, \sum_{\substack{\forall p_{tc} \in \\ P_{r,tc}}} d_{t'',p_{tc}}^{t'} \!\! - M \, y_{t',tc}^t - \\ \!\!\!\!\!\! M \left(1 - \!\! \sum_{\substack{\forall r \in R^t: \\ tc \in TC_r}} \!\! x_{r}^{t} \, ebs_{r,tc}^t \right) - M \left(1 -  \!\! \sum_{\substack{\forall r \in R^{t'}: \\ tc \in TC_r}}  \!\! x_{r}^{t'}\right) \label{eq:submodel delay4}
\end{multline}

\begin{multline}
        sR_{tc}^{t} \geq eR_{tc}^{t'} + \!\! \sum_{\substack{\forall t'' \in T: \\ t'' \neq t}} \, \sum_{\substack{\forall r \in R^{t'}: \\ tc \in TC_r}} \!\!\!\!\!\!\!\!\!\! \sum_{\:\:\:\:\: \substack{\forall tc^{\star} \in TC_r: \\ \:\:\:\:\: s_{r,tc} = ref_{r,tc^{\star}} }} \!\!\!\!\!\!\!\!\!\!\!\!\! d_{t'',s_{tc}}^{t'} \!\! - M \, y_{t',tc}^t - \\ M \left(1 - \!\! \sum_{\substack{\forall r \in R^t: \\ tc \in TC_r}} \!\! x_{r}^{t} \, ebs_{r,tc}^t \right) - M \left(1 -  \!\! \sum_{\substack{\forall r \in R^{t'}: \\ tc \in TC_r}}  \!\! x_{r}^{t'}\right) \label{eq:submodel delay5}
\end{multline}
\begin{algorithm}
    \caption{Precedence constraints}
    \label{alg:conflict constraints}
    \begin{algorithmic}[1] 
        \For{$\forall t \in T$}
            \For{$\forall tc \in TC^t$}
                \For{$\forall t' \in T: t' \neq t \land tc \in TC^{t'}$}
                    \For{$r \in R^t: tc \in TC_r$}
                        \If{$sw_{r,tc}^c = 0 \land sw_{r,p_{r,tc}}^o = 0$}
                            \State $y_{t',tc}^t = y_{t',p_{r,tc}}^t$ \label{constr: prec3}
                        \ElsIf{$sw_{r,tc}^c = 1$}
                            \State $y_{t',tc}^t - y_{t',p_{r,tc}}^t - \!\!\!\! \sum \limits_{\substack{\forall r' \in R^t: \\ p_{r,tc} \in TC_{r'}}} \!\!\!\! x_r^t \geq -1$ \label{constr: prec4}
                            \vspace{1mm}
                            \State $y_{t',tc}^t - y_{t',p_{r,tc}}^t + \!\!\!\! \sum \limits_{\substack{\forall r' \in R^{t'}: \\ p_{r,tc} \in TC_{r'}}} \!\!\!\! x_{r'}^{t'} \leq 1$ \label{constr: prec5}
                        \Else
                            \State $y_{t',tc}^t - y_{t',p_{r,tc}}^t - \! \sum \limits_{\substack{\forall r' \in R^{t'}: \\ tc \in TC_{r'}}} \!\! x_{r'}^{t'} \geq -1$ \label{constr: prec6}
                            \vspace{1mm}
                            \State $y_{t',tc}^t - y_{t',p_{r,tc}}^t + \! \sum \limits_{\substack{\forall r' \in R^t: \\ tc \in TC_{r'}}} \! x_r^t \leq 1$ \label{constr: prec7}
                        \EndIf
                    \EndFor
                    
                    \For{$t'' \in T: t'' \neq t, t'' \neq t'\land tc \in TC^{t''}$}
                        \State $y_{t',tc}^t + y_{t'',tc}^{t'} \leq 1 + y_{t'',tc}^t$ \label{constr: prec8}
                        \State $y_{t',tc}^t + y_{t'',tc}^{t'} \geq y_{t'',tc}^t$ \label{constr: prec9}
                    \EndFor
                    \vspace{1mm}
                    \State $y_{t',tc}^t \geq \sum \limits_{\substack{\forall r \in R^{t'}: \\ tc \in TC_r}} x_r^{t'} - \!\! \sum \limits_{\substack{\forall r \in R^{t}: \\ tc \in TC_r}} x_r^{t}$ \label{constr: prec10}
                    \State $y_{t,tc}^{t'} \geq \sum \limits_{\substack{\forall r \in R^{t}: \\ tc \in TC_r}} x_r^{t} - \!\! \sum \limits_{\substack{\forall r \in R^{t'}: \\ tc \in TC_r}} x_r^{t'}$ \label{constr: prec11}
                \EndFor
            \EndFor
        \EndFor
    \end{algorithmic}
\end{algorithm}

\subsection{Optimality Condition Analysis} \label{section:optimality condition}
\begin{figure*}[t]
    \centering
    \includegraphics[width=0.99\linewidth]{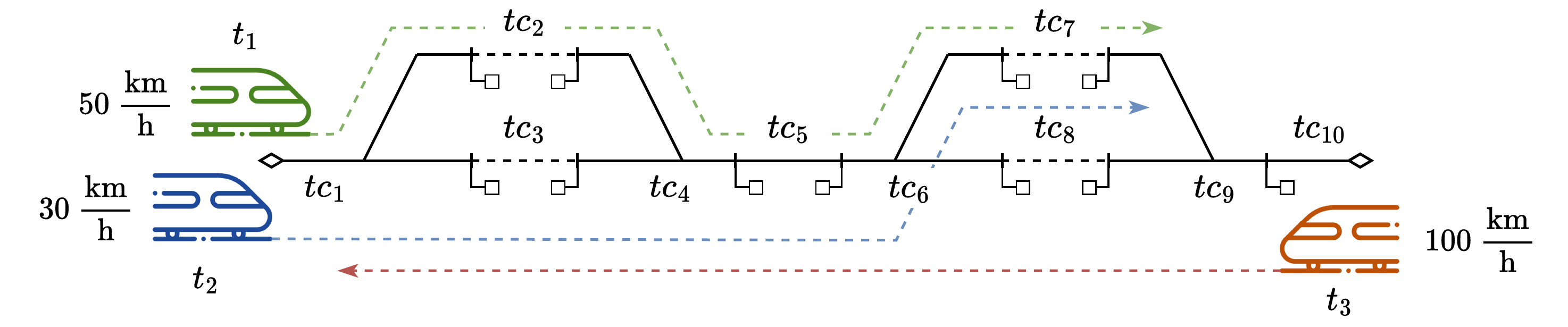}
    \caption{Example scenario resulting in sub-optimal solution}
    \label{fig:suboptimal example}
\end{figure*}

The conditions in which the model detailed in the previous section provides a sub-optimal solution is investigated and given explicitly in this section. Since the sub-optimal model tackles the rerouting and reordering problem as well, the optimality condition of the reordering is discussed first with fixed route selection. An example scenario involving three trains on a network modelling two connected stations is illustrated in Fig.~\ref{fig:suboptimal example} when the solution of the model is sub-optimal in terms of train precedence. The platforms represented by dashed lines are 400 m in length with 100 and 40 km/h maximum velocity on the lower and upper branches, respectively. The track circuits, including the switches, have 200 m length along each route, while $tc5$ connecting the stations is 1200 m long, and the maximum speed allowed on them is uniformly set to 100 km/h. The trains with $w^t = 1$ identical priorities would travel with constant speed in accordance with Fig.~\ref{fig:suboptimal example}, but $t_1$ that is scheduled to dwell two additional minutes on platform $tc_7$. According to the sub-optimal model, train $t_3$ highlighted with red should wait for both $t_1$ and $t_2$ at track circuit $tc_5$, and $t_2$ is delayed at $tc_3$ due to the two-minute dwell of $t_1$. Therefore, the objective value $D_{\Sigma}(X_1)$ of the solution $X_1$ is:
\begin{equation}
    D_{\Sigma}(X_1) = d_{t_1, tc_8}^{t_3} + d_{t_2, tc_8}^{t_3} + d_{t_1, tc_3}^{t_2} \, .
\end{equation}
However, since the $d_{t_1, tc_3}^{t_2}$ delay assigned to $t_2$ is propagated to $t_3$ as well, the $\widetilde{D}_{\Sigma}(X_1)$ total delay of the trains would be:
\begin{equation}
    \widetilde{D}_{\Sigma}(X_1) = 2 \, d_{t_1, tc_3}^{t_2} + d_{t_2, tc_8}^{t_3} .
\end{equation}
In the optimal solution $X_2$, the green train $t_1$ is overtaken by $t_2$ at track circuit $tc_4$, also giving the precedence to $t_3$ with $\widetilde{D}_{\Sigma}(X_2) < \widetilde{D}_{\Sigma}(X_1)$ global objective value:
\begin{equation}
    \widetilde{D}_{\Sigma}(X_2) = 2 \, d_{t_2, tc_8}^{t_3} + d_{t_3, tc_2}^{t_1} .
\end{equation}
Since however, the objective value $D_{\Sigma}(X_2)$ of the global optimum $X_2$ according to the sub-optimal model as:
\begin{equation}
    D_{\Sigma}(X_2) = d_{t_2, tc_2}^{t_1} + d_{t_2, tc_8}^{t_3} + d_{t_3, tc_2}^{t_1} > D_{\Sigma}(X_1) \, ,
\end{equation}
the solution $X_1$ is sub-optimal in terms of total delay. 

In general, the solution $X$ of the reordering or the rerouting problems neglecting the delay propagation is sub-optimal compared to the original model without overlaps if:
\begin{align}
   \begin{split}
   \exists \, {y_{t',tc}^{t}}' \in X' \quad &\forall \{t,t'\} \in T, \forall tc \in TC^t \cap TC^{t'} \, \lor \\
   \exists \, {x_{r}^{t}}' \in X' \quad &\forall t \in T, \forall r \in R
   \end{split} \label{eq: alter solution}
\end{align}
satisfying the constraints in Alg.~\ref{alg:capacity constraints} and the resulting ${d_{t',tc}^{t}}'$ delays according to \eqref{eq:submodel delay1} and \eqref{eq:submodel delay2} meet the following criteria:
\begin{align}
    f_s \left({d_{t',tc}^{t}}', \, {\delta^t}' \, \rvert \, X' \right) &> f_s \Bigl( d_{t',tc}^{t}, \, \delta^t \, \rvert \, X \Bigr) , \label{eq:suboptimal1} \\ 
    f_s \left( \widetilde{d}_{t',tc_r^{\star}}^{t}{}', \, {\delta^t}' \, \rvert \, X' \right) &\leq f_s \left(\widetilde{d}_{t',tc_r^{\star}}^{t}, \, \delta^t \, \rvert \, X\right) . \label{eq:suboptimal2}
\end{align}
where $\widetilde{d}_{t',tc_r^{\star}}^{t}$ denotes the delay assigned to train $t$ at track circuit $tc_r^{\star} = p_{r, ref_{r,tc}}$ due to the conflict with $t'$ at $tc$, including the delay propagation as:
\begin{equation}
    \sum_{\substack{\forall r \in R^t: \\ tc \in TC^r}} \widetilde{d}_{t',tc_r^{\star}}^{t} = eR_{tc}^{t'} - sR_{tc}^{t} + \, ^{\uparrow}\!\Delta_{t',tc}^{t} - \, ^{\downarrow}\!\Delta_{t',tc}^{t}\,. \label{eq:mod delay}
\end{equation}
Constraint \eqref{eq:mod delay} considers the signaling system by delaying train $t$ at track circuit $p_{r, ref_{r,tc}}$ that terminates the block section preceding $tc$ alogn route $r$. The propagation variables $^{\uparrow}\!\Delta_{t',tc}^{t}$ and $^{\downarrow}\!\Delta_{t',tc}^{t}$ in \eqref{eq:mod delay} increasing and decreasing the $d_{t',tc}^{t}$ delay of train $t$ at track circuit $tc$ with respect to $t'$ is quantified by the following equations:
\begin{align}
    &^{\uparrow}\!\Delta_{t',tc}^{t} = \!\! \sum_{\substack{\forall t'' \in T : t'' \neq t', \\ tc \in TC^{t''}}} \,\, \sum_{\substack{\forall r \in R^{t'}: \\ tc \in TC_r}} \,\, \sum_{\substack{\forall p_{tc} \in P_{r,tc} \setminus t'' = t, \\ p_{tc} \in p_{r, ref_{r,tc}}}} \!\!\!\!\! d_{t'',p_{tc}}^{t'} \, , \\
    &^{\downarrow}\!\Delta_{t',tc}^{t} = \!\! \sum_{\substack{\forall t'' \in T : t'' \neq t, \\ tc \in TC^{t''}}} \,\, \sum_{\substack{\forall r \in R^{t}: \\ tc \in TC_r}} \,\, \sum_{\substack{\forall p_{tc} \in P_{r,tc} \setminus t'' = t', \\ p_{tc} \in p_{r, ref_{r,tc}}}} \!\!\!\!\! d_{t'',p_{tc}}^{t} \, ,
\end{align}
so that substituting $d_{r, tc}^t = \sum_{t' \in T: t' \neq t} \widetilde{d}_{t',tc}^{t} \, x_{r}^t$ in Alg.~\ref{alg:scheduling constraints} and Alg.~\ref{alg:capacity constraints} results in a feasible solution. Since the subtrahend $^{\downarrow}\!\Delta_{t',tc}^{t}$ in \eqref{eq:mod delay} decreases the objective value, false delays may occur, leading to an infeasible solution according to the original model. This problem is tackled by the global optimum model detailed in the following section.

\subsection{Global Optimum Model} \label{section:global optimum model}

\begin{algorithm}
    \caption{Global optimum delays}
    \label{alg:global constraints}
    \begin{algorithmic}[1] 
        \For{$\forall t \in T$}
            \For{$\forall tc \in TC^t$}
                \For{$\forall t' \in T: t' \neq t \land tc \in TC^{t'}$}
                    \State $\widetilde{d}_{t',tc}^{t} = M  \!\!\!\!\!\!\!\! \sum \limits_{\substack{r \in R^{t}: \\ bs_{r,tc} \neq bs_{r,s_{r,tc}}}} \!\!\!\!\!\!\!\! x_{r}^{t'} $ \label{constr: glob1}
                    \vspace{1mm}
                    \State $\widetilde{d}_{t,tc}^{t'} = M  \!\!\!\!\!\!\!\! \sum \limits_{\substack{r \in R^{t'}: \\ bs_{r,tc} \neq bs_{r,s_{r,tc}}}} \!\!\!\!\!\!\!\! x_{r}^{t} $ \label{constr: glob2}
                    
                    \vspace{1mm}
                    
                    \State $\sum \limits_{\substack{\forall r \in R^t: \\ tc \in TC^r}} \widetilde{d}_{t',p_{r, ref_{r,tc}}}^{t} \leq M(1 - y_{t',tc}^{t}) \: + $ \label{constr: glob3}
                    \vspace{-3mm}
                    \Statex \hspace{40mm} $+ M \left(1 - \sum \limits_{{\substack{r \in R^{t}: \\ tc \in TC^r}}} \!\!\! x_{r}^{t} \right)$ 
                    \vspace{3mm}
                    \State $\sum \limits_{\substack{\forall r \in R^t: \\ tc \in TC^r}} \widetilde{d}_{t,p_{r, ref_{r,tc}}}^{t'} \leq M(1 - y_{t,tc}^{t'}) \: + $ \label{constr: glob4}
                    \vspace{-3mm}
                    \Statex \hspace{40mm} $+ M \left(1 - \sum \limits_{{\substack{r \in R^{t'}: \\ tc \in TC^r}}} \!\!\! x_{r}^{t'} \right)$
                    
                    \vspace{2mm}
                    
                    \State $\sum \limits_{\substack{\forall r \in R^t: \\ tc \in TC^r}} \widetilde{d}_{t',tc}^{t} \leq 
                    M \sum \limits_{\substack{\forall r \in R^t: \\ tc \in TC^r}} x_{r}^{t}$ \label{constr: glob5}
                    \vspace{1mm}
                    \State $\sum \limits_{\substack{\forall r \in R^{t'}: \\ tc \in TC^r}} \widetilde{d}_{t,tc}^{t'} \leq 
                    M \sum \limits_{\substack{\forall r \in R^{t'}: \\ tc \in TC^r}} x_{r}^{t'}$ \label{constr: glob6}
                    \vspace{2mm} 
                    \For{$\forall r \in R^t: tc \in TC_r$}
                    \vspace{1mm} 
                        \State $\widetilde{d}_{t', tc}^{t} \leq 
                        M \!\!\!\!\!\!\! \sum \limits_{\substack{\forall r' \in R^{t'}: \\ bs_{r,s_{r,tc}} \cap r' \neq 0 }} \!\!\!\!\!\!\! x_{r'}^{t'}$ \label{constr: glob7}
                    \EndFor
                    \For{$\forall r \in R^{t'}: tc \in TC_r$}
                    \vspace{1mm} 
                        \State $\widetilde{d}_{t, tc}^{t'} \leq 
                        M \!\!\!\!\!\!\! \sum \limits_{\substack{\forall r' \in R^{t}: \\ bs_{r,s_{r,tc}} \cap r' \neq 0 }} \!\!\!\!\!\!\! x_{r'}^{t}$ \label{constr: glob8}
                    \EndFor
                    
                    \For{$t'' \in T: t'' \neq t, t'' \neq t'\land tc \in TC^{t''}$}
                        
                        \State\vspace*{-\baselineskip}
                        \begin{fleqn}[\dimexpr(\leftmargini-\labelsep)*4]
                        \setlength\belowdisplayskip{0pt}
                        \begin{equation*}
                            \begin{multlined}[c]
                                \textstyle \sum \limits_{\substack{\forall r \in R^{t}: \\ tc \in TC^r}} \!\! \widetilde{d}_{t',p_{r, ref_{r,tc}}}^{t} \!\! \leq M \!\! \left( 2 - y_{t'',tc}^{t'} - y_{t,tc}^{t''} \right) \!\!\!\!\!\!\!\! \vspace{-2mm} \\
                                + M \left(2 - \textstyle \sum \limits_{\substack{\forall r \in R^{t'}: \\ tc \in TC^r}} x_{r}^{t'} - \textstyle \sum \limits_{\substack{\forall r \in R^{t''}: \\ tc \in TC^r}} x_{r}^{t''} \right)
                            \end{multlined}
                        \end{equation*}
                        \end{fleqn} \label{constr: glob9}
                        
                        \vspace{2mm}
                        
                        \State\vspace*{-\baselineskip}
                        \begin{fleqn}[\dimexpr(\leftmargini-\labelsep)*4]
                        \setlength\belowdisplayskip{0pt}
                        \begin{equation*}
                            \begin{multlined}[c]
                                \textstyle \sum \limits_{\substack{\forall r \in R^{t'}: \\ tc \in TC^r}} \!\! \widetilde{d}_{t,p_{r, ref_{r,tc}}}^{t'} \!\! \leq M \!\! \left(2 - y_{t'',tc}^{t} - y_{t',tc}^{t''} \right) \!\!\!\!\!\!\!\! \vspace{-2mm} \\
                                + M \left( 2 - \textstyle \sum \limits_{\substack{\forall r \in R^{t}: \\ tc \in TC^r}} x_{r}^{t} - \textstyle \sum \limits_{\substack{\forall r \in R^{t''}: \\ tc \in TC^r}} x_{r}^{t''} \right)
                            \end{multlined}
                        \end{equation*}
                        \end{fleqn} \label{constr: glob10}
                    \EndFor
                \EndFor   
            \EndFor
        \EndFor
    \end{algorithmic}
\end{algorithm}

The proposed global optimum model extends the formulation of $d_{t',tc}^{t}$ in the sub-optimal model involving the delay propagation according to the optimality condition analysis as:  
\begin{multline}
     \sum_{\substack{\forall r \in R^t: \\ tc \in TC^r}} \widetilde{d}_{t',p_{r, ref_{r,tc}}}^{t} \geq eR_{tc}^{t'} - sR_{tc}^{t} + \, ^{\uparrow}\!\Delta_{t',tc}^{t} - \, ^{\downarrow}\!\Delta_{t',tc}^{t} - \\ - M \, y_{t',tc}^t - M \left(1 - \sum_{\substack{\forall r \in R^{t'}:\\ tc \in TC_r}} x_{r}^{t'}\right) ,
\end{multline}
where the last two terms neglect to delay train $t$ at track circuit $tc$ where $t$ has priority over $t'$ or it does not use $tc$. However, false values decreasing the objective value may occur due to the $^{\downarrow}\!\Delta_{t',tc}^{t}$ total delay accumulated by train $t$ to track circuit $tc$. The false delays result in an infeasible solution is prevented by the constraints in Alg.~\ref{alg:global constraints}. Besides the delay propagation, the global optimum model considers the signaling system of the control area concerning the block sections. Therefore, the potential conflict of trains $t$ and $t'$ at track circuit $tc$ is resolved at the track circuit $p_{r, ref_{r,tc}}$ preceding the corresponding block section $bs_{r,tc}$, and the rescheduling at track circuits lacking a signal is prevented by the constraints in lines \ref{constr: glob1} - \ref{constr: glob2} of Alg.~\ref{alg:global constraints}. The inequalities in lines \ref{constr: glob3} - \ref{constr: glob4} impose that the delay assigned to train $t$ due to the conflict with $t'$ at track circuit $tc$ can only be nonzero if $t'$ is reported to use $tc$ before $t$. According to lines \ref{constr: glob5} - \ref{constr: glob6}, the trains cannot be delayed on track circuits they do not use. A conflict between two trains may arise if a block section is reserved by both of them. Hence, the false delay assigned to train $t$ at track circuit $tc$ due to $t'$ is avoided by the constraints in lines \ref{constr: glob7} and \ref{constr: glob8} when the block section $bs_{r,s_{r,tc}}$ of the subsequent track circuit $s_{r,tc}$ does not intersect the $r'$ route of $t'$. Finally, according to expression in lines \ref{constr: glob9} - \ref{constr: glob10}, if train $t$ has to wait for multiple trains at a given track circuit, the delay assigned to it is computed from the schedule of the last train to pass. Although this model provides a global optimum for the original problem in \cite{pellegrini_2012,pellegrini2015recife}, it is not guaranteed for the proposed model of the rtRTMP extended with the overlap constraints. The global optimum model can result in a sub-optimal solution $X$ when a feasible alternative solution $X'$ in accordance with \eqref{eq: alter solution} so that:
\begin{align}
    f_s \left( \widetilde{d}_{t',tc_r^{\star}}^{t}{}', \, \delta^t \, \rvert \, X' \right) &\geq f_s \left(\widetilde{d}_{t',tc_r^{\star}}^{t}, \, \delta^t \, \rvert \, X\right) , \label{eq:suboptimal3} \\ 
    f_s \left( \widetilde{d}o_{t',tc_r^{\star}}^{t}{}', \, {\delta^t}' \, \rvert \, X' \right) &\leq f_s \left(\widetilde{d}o_{t',tc_r^{\star}}^{t}, \, {\delta^t}' \, \rvert \, X\right) , \label{eq:suboptimal4}
\end{align}
where $\widetilde{d}o_{t',tc_r^{\star}}^{t}$ denotes the delay that has to be assigned to train $t$ to resolve the conflict with $t'$ at track circuit $tc$ considering the overlaps according to Alg. \ref{alg:overlap constraints} as:
\begin{multline}
      \widetilde{d}o_{t',tc_r^{\star}}^{t} = \widetilde{d}_{t',tc_r^{\star}}^{t} \!+\! y_{t',tc}^{t} \Bigl( \Bigl(1 - c_{t',tc}^{t} \Bigr) \! \left(eR_{tc}^{t'} - sO_{tc}^{t} \right) \, + \\ c_{t',tc}^{t} \left( eO_{tc}^{t'} - sR_{tc}^{t} \! \right) \Bigr) \!+\! \left(\! sO_{tc}^{t'} - eO_{tc}^{t} \! \right) yO_{t',tc}^{t}.
\end{multline}
The optimality of the proposed global optimum model concerning the extended rtRTMP is investigated experimentally in the following section.

\section{Experimental Results} \label{section:results}
\begin{figure*}
    \centering
    \subfloat[Network A] {\includegraphics[width=0.99\linewidth]{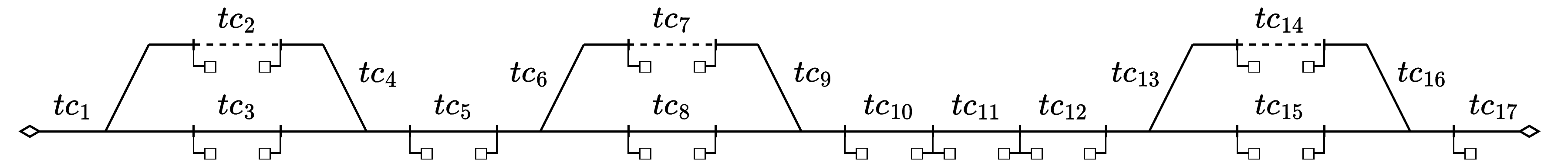} \label{fig:infra A}}
    
    \hfil
    
    \subfloat[Network B] {\includegraphics[width=0.99\linewidth]{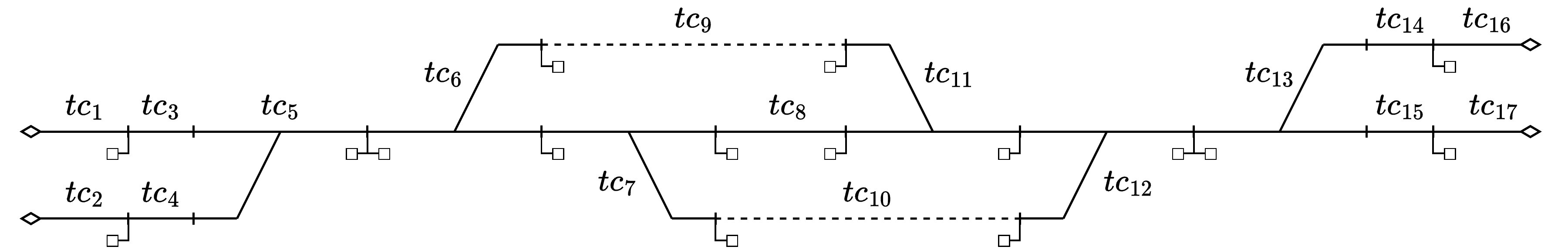} \label{fig:infra B}}
    
    \caption{Infrastructure model of the rail networks used for evaluation}
    \label{fig: infras}
\end{figure*}
The proposed models and the multi-stage optimization workflow is evaluated in two different rail networks and various scenarios. The infrastructure model of the two control areas illustrated in Fig.~\ref{fig: infras} is constructed manually in a Matlab GUI implemented for this purpose. The dashed lines in Fig.~\ref{fig: infras} represent the platforms, and the signals are located at the boundaries of the track circuits. The network in Fig.~\ref{fig:infra A}, consisting of 17 track circuits and eight routes, representing three stations connected by single rail sections. The infrastructure in Fig.~\ref{fig:infra B} also includes 17 track circuits, forming 12 routes between the four entry and exit locations. The details of the track circuits, including their total length and speed limit, are summarized in Table~\ref{infra params}. Two hundred scenarios are generated randomly in both networks, one hundred involving three trains and the other one hundred with four trains. The parameters defining the initial state and properties of the trains are uniformly distributed between the limits given in Table~\ref{train params}. The trains can enter the control area at either of the entry locations, or may have already been in the network at the beginning of the optimization. The trains that entered the control area before, referred to as initial trains, are drawn with a given $P_{in}$ discrete probability. The location of the initial trains is selected from the track circuits with even probability. The trains imposed to wait with probability $P_{dw}$ must dwell some uniformly distributed random time at either of the platforms. The randomly generated scenarios must be feasible according to the original model defined in Section \ref{section:Model Formulation} extended with the overlap constraints. The proposed models are evaluated in terms of the objective value and runtime of the optimization. The performance of the proposed models is given relatively to the original models averaging the results of the $N$ different scenarios with the following metrics:
\begin{align}
    P_{obj} &= \sum_{i = 1}^{N} \left(1 + \dfrac{f(X_i) - f(\hat{X}_i)}{\sum_{i = 1}^{N} f(\hat{X}_i)} \right) \cdot 100 \, , \label{eq:obj value metric} \\
    P_{run} &= \dfrac{1}{N} \sum_{i = 1}^{N} \dfrac{rt_i}{\hat{rt}_i} \cdot 100 \, , \label{eq:runtime metric}
\end{align}
where $f(X_i)$, $f(\hat{X}_i)$, $rt_i$ and $\hat{rt}_i$ denotes the objective value and runtime of the proposed and reference model in the $i$-th scenario. The runtime performance is given simply by the average ratio of the proposed and reference model in \eqref{eq:runtime metric}. However, since the objective value can be zero, the metric in \eqref{eq:obj value metric} is defined based on the difference between the proposed and reference model and the average objective value of the reference model. The objective value is computed according to the reference model as: 
\begin{equation}
    f(X) = \sum_{\forall t \in T} w^t D^t.
\end{equation}
Both proposed models solve either the rerouting and reordering sub-problems of the rtRTMP, determining the local delays of the trains based on the original schedule and resulting precedence indicators as:
\begin{multline}
    d_{r, tc^{\star}}^{t} = 
    \max \limits_{\substack{\forall t' \in T: t' \neq t \\ \forall tc \in bs_{r,s_{r,tc}}}} \!\! \Biggl(
    eR_{tc}^{t'} - sR_{tc}^{t} -
    \sum \limits_{\substack{\forall p_{tc^{\star}} \in \\ P_{r,tc^{\star}}}} \! d_{r,p_{tc^{\star}}}^{t} + \\
    \sum \limits_{\substack{r' \in R^{t'} : \\ tc \in r'}} \, \sum \limits_{\substack{\forall p_{tc}' \in \\ P_{r',tc}}} \!\! d_{r',p_{tc}'}^{t'} \Biggr) x_{r}^t \, y_{t,tc}^{t'} \, sig_{r,tc^{\star}} \, . \label{eq: sched from prec}
\end{multline}
In Eq. \eqref{eq: sched from prec}, $sig_{r,tc^{\star}} \in \{0,1\}$ indicates if $tc$ terminates its $bs_{r,tc}$ block section with a signal along route $r$ as: 
\begin{equation}
    sig_{r,tc^{\star}} = \begin{cases}
        1 & \text{if\ } bs_{r,tc} \neq bs_{r,s_{r, tc}} \\
        0 & \text{otherwise}
    \end{cases} \, .
\end{equation}
Therefore, the first stage of the optimization is compared directly with the original model neglecting the overlaps. Despite the constraints aimed at this, the sub-optimal model may provide an infeasible solution in some cases due to the simplified formulation of the delay propagation in \eqref{eq:submodel delay4} and \eqref{eq:submodel delay5}. In this case, the model is evaluated as if it provides the solution of the reference model, increasing the response time by the runtime of the sub-optimal model. The reordering and rescheduling are not bijective sub-problems, considering the overlaps. The second stage of the optimization is intended to solve the rescheduling sub-problem of the rtRTMP extended with overlap constraints in Alg. \ref{alg:overlap constraints} using the solution of the first stage provided by the proposed models. Therefore, the reference model of the second stage is extended with the variables and constraints related to the overlaps without using the preliminary solution. Since multiple information can be extracted from the solution of the first stage, the second stage is assessed in two ways. First, only the route selection according to the first stage is applied, and then the precedence of the trains is considered too. The overlap constraints may not be feasible with the rerouting or reordering reported by the proposed models in the first stage. This is assessed similarly to the infeasible solution of the sub-optimal model in the first stage. The evaluation is performed on a Lenovo ThinkCentre PC with an Intel Core i7-10700 2.9 GHz processor and 16 GB memory. The optimization problems are solved in Matlab environment by applying branch and bound algorithm. Since it requires the initial condition to be complete and feasible, the solution to the rerouting and reordering problems form only hard constraints in the evaluation of the second stage. Moreover, the runtime decrease of using an initial condition depends highly on the built-in heuristics of the MILP solver.

The evaluation results are given separately for the two networks in Table~\ref{eval results A} and Table~\ref{eval results B}, including the relative objective value and runtime performance of the proposed models in the two stages of the optimization workflow. The computational complexity of the optimization is decreased significantly by both of the proposed models. In most circumstances, the runtime of the sub-optimal model is around 50\% of the reference model in the first stage of the optimization. Even in the most complex scenarios, the proposed model is faster by 39\% on average. Despite the significant runtime improvement, the maximum overhead of the model simplification regarding the objective value is 25\%. The 52\textendash43\% runtime decrease in the second optimization stage is even more significant using the preliminary solution of the first stage provided by the sub-optimal model. However, since constraining the precedence between the trains besides the route selection leads more frequently to an infeasible solution considering the overlaps, sometimes it even increases the runtime. Although the global optimum model does not decrease the response time of the first stage as much as the sub-optimal one, it always provides the same solution as the reference model 27\textendash46\% faster. Even the objective value in the second stage does not deteriorate more than 17\%. Despite the complexity of the global optimum model, it is slightly more efficient than the sub-optimal one with respect to the runtime due to the faster convergence and lower probability of an infeasible solution. 

\begin{table}[t]
\caption{Evaluation results in infrastructure model of network A}
\label{eval results A}
\centering
\setlength\extrarowheight{2.5pt}
\begin{tabular}{@{}lcccccc@{}} 
\toprule
 & \multicolumn{3}{c}{Sub-optimal} & \multicolumn{3}{c}{Global optimum} \\
 \cmidrule(lr){2-4} \cmidrule(ll){5-7}
 & \multirow{2}{*}{1\textsuperscript{st}} & \multicolumn{2}{c}{2\textsuperscript{nd}} & \multirow{2}{*}{1\textsuperscript{st}} & \multicolumn{2}{c}{2\textsuperscript{nd}}\\
  \cmidrule(lr){3-4} \cmidrule(ll){6-7}
 & & Route & Precedence & & Route & Precedence\\
\midrule
\multicolumn{7}{c}{3 Trains}\\
\midrule
$P_{obj}$ [\%] & $125$ & $109$ & $120$  & $100$ & $102$ & $102$ \\
$P_{run}$ [\%] & $48$ & $31$ & $28$ & $54$ & $30$ & $27$ \\
\midrule
\multicolumn{7}{c}{4 Trains}\\
\midrule
$P_{obj}$ [\%] & $118$ & $108$ & $113$ & $100$ & $101$ & $102$ \\
$P_{run}$ [\%] & $50$ & $27$ & $35$ & $65$ & $27$ & $25$ \\
\bottomrule
\end{tabular}
\end{table}

\begin{table}[t]
\caption{Evaluation results in infrastructure model of network B}
\label{eval results B}
\centering
\setlength\extrarowheight{2.5pt}
\begin{tabular}{@{}lcccccc@{}} 
\toprule
 & \multicolumn{3}{c}{Sub-optimal} & \multicolumn{3}{c}{Global optimum} \\
 \cmidrule(lr){2-4} \cmidrule(ll){5-7}
 & \multirow{2}{*}{1\textsuperscript{st}} & \multicolumn{2}{c}{2\textsuperscript{nd}} & \multirow{2}{*}{1\textsuperscript{st}} & \multicolumn{2}{c}{2\textsuperscript{nd}}\\
  \cmidrule(lr){3-4} \cmidrule(ll){6-7}
 & & Route & Precedence & & Route & Precedence\\
\midrule
\multicolumn{7}{c}{3 Trains}\\
\midrule
$P_{obj}$ [\%] & $118$ & $115$ & $126$  & $100$ & $112$ & $117$ \\
$P_{run}$ [\%] & $50$ & $48$ & $48$ & $55$ & $47$ & $45$ \\
\midrule
\multicolumn{7}{c}{4 Trains}\\
\midrule
$P_{obj}$ [\%] & $125$ & $115$ & $124$ & $100$ & $108$ & $111$ \\
$P_{run}$ [\%] & $61$ & $32$ & $27$ & $73$ & $30$ & $22$ \\
\bottomrule
\end{tabular}
\vspace{-4mm}
\end{table}

\begin{figure}
    \centering
    \includegraphics[width=0.99\linewidth]{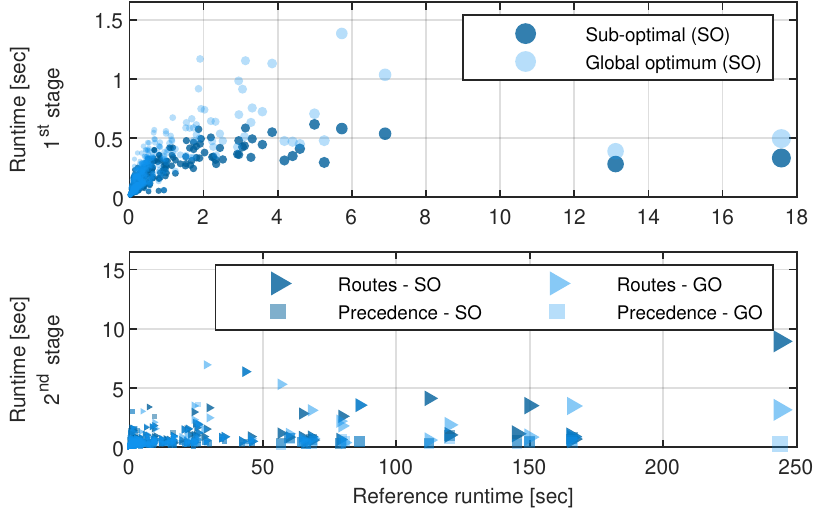}
    \caption{Runtime comparison between the proposed and reference models}
    \label{fig:runtime comparison}
\end{figure}

\begin{figure}
    \centering
    \includegraphics[width=0.99\linewidth]{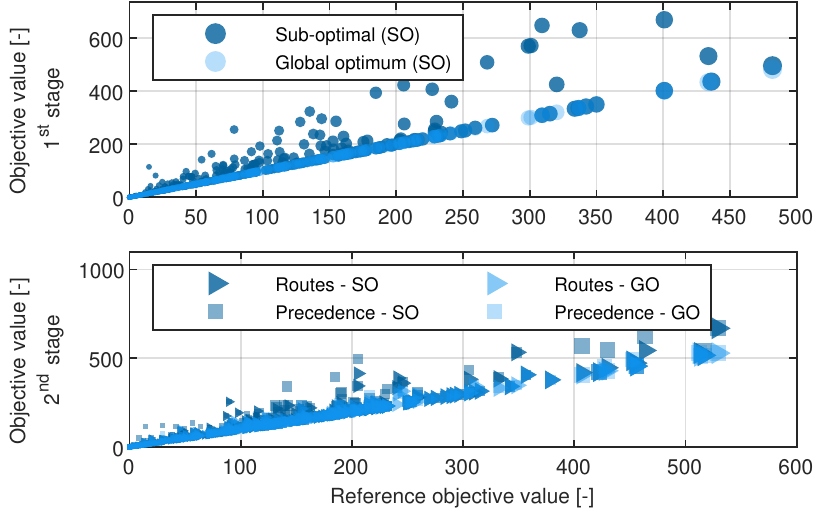}
    \caption{Objective value comparison between the proposed and reference models}
    \label{fig:obj value comparison}
     \vspace{-5mm}
\end{figure}

The average runtime decrease varies between 27\% and 78\% in the different circumstances. However, the runtime complexity trend of the proposed algorithms illustrated in Fig.~\ref{fig:runtime comparison} is even more favourable. The more complex the scenario is, the more significant the impact of the proposed models solving the problem with the highest 15-second runtime more than 35 times faster than the original model. At the same time, the runtime peak of the sub-optimal and global optimum models are less than 1.5 and 0.75 seconds, respectively, neglecting the overlap constraints. In the second stage of the optimization, the 244-second maximum response time of the reference model, is decreased by 96\textendash99.9\% by using the preliminary route selection and precedence indicators reported by the global optimum model. The objective value of the proposed models compared to the reference depicted in Fig~\ref{fig:obj value comparison} does not show the same tendency as the runtime complexity. In most cases, the solution of the sub-optimal model differs slightly from the reference model. Although the train delays weighted by their priorities are twice as higher as the reference in some scenarios neglecting the overlap constraints, the relative performance of the proposed model does not change with the objective value of the reference model. Besides granting the solution that minimizes the delay propagation without overlaps, the global optimum model also provides good input for the second stage of the optimization, resulting in similar objective values to the reference. 

\section{Conclusion}
The reinterpreted delays resolving the local conflict between two trains allow the formulation of the proposed models based on the original timetable without any feedback. The simplified sub-optimal model provides a fast solution for the original problem without the overlap constraints. Despite neglecting the delay propagation, the performance of the optimization does not degrade significantly on average. The global optimum model reduces the average runtime almost as much as the sub-optimal one, ensuring the best solution in terms of the optimization objective. The MILP model of the real-time railway traffic management problem extended with the safety-critic regulation of the overlaps may not be able to resolve the conflicts within the desired response time. While the proposed multi-stage optimization workflow, decomposing the rerouting, reordering, and rescheduling problems, is characterized by a relatively low increase of the objective value, it significantly reduces the complexity of the problem by using the preliminary solution of the proposed models. Furthermore, the higher runtime of the reference model, the more efficient decrease can be achieved by the proposed models, resulting in a much more favourable complexity class. Therefore, the multi-stage optimization workflow with the sub-optimal and global-optimum models can sufficiently solve the extended real-time railway problem. 

\section*{Appendix}

\begin{table}[H]
\caption{Train parameter settings}
\label{train params}
\centering
\setlength\extrarowheight{2.5pt}
\begin{tabular}{@{}llcc@{}} 
\toprule
Parameter & Notation & Lower limit & Upper limit\\
\midrule
Desired velocity & $v^t$ [km/h] & $50$ & $160$ \\
Length of vehicle assembly & $L_v^t$ [m] & $100$ & $300$ \\
Train priority & $w^t$ [\%] & $0.1$ & $100$ \\
Scheduled entry time & $init^t$ [sec] & $0$ & $50$ \\
Waiting time & $dw^t$ [sec] & $0$ & $180$ \\
Waiting probability & $P_dw$ [\%] & \multicolumn{2}{c}{$25$} \\
Initial train probability & $P_{in}$ [\%] & \multicolumn{2}{c}{$20$} \\
\bottomrule
\end{tabular}
\end{table}

\begin{table}[H]
\caption{Parameters of the infrastructure models}
\label{infra params}
\centering
\setlength\extrarowheight{2.5pt}
\begin{tabular}{@{}llccccc@{}} 
\toprule
 \multicolumn{2}{l}{Track} & \multicolumn{2}{c}{Network A} & & \multicolumn{2}{c}{Network B}\\
\cline{3-4}
\cline{6-7}
\multicolumn{2}{l}{circuit} & Vel. [km/h] & Length [m] & & Vel. [km/h] & Length [m]\\
\midrule
\multirow{2}{*}{$tc_1$} & $\uparrow$ & $100$ & $100$ & & $40$ & $200$ \\
& $\downarrow$ & $-$ & $-$ & & $100$ & $200$ \\
\multirow{2}{*}{$tc_2$} & $\uparrow$ & $-$ & $-$ & & $40$ & $400$ \\
& $\downarrow$ & $100$ & $100$ & & $-$ & $-$ \\
\multirow{2}{*}{$tc_3$} & $\uparrow$ & $100$ & $50$ & & $-$ & $-$ \\
& $\downarrow$ & $-$ & $-$ & & $100$ & $400$ \\
\multirow{2}{*}{$tc_4$} & $\uparrow$ & $-$ & $-$ & & $40$ & $200$ \\
& $\downarrow$ & $100$ & $50$ & & $100$ & $200$ \\
\multirow{2}{*}{$tc_5$} & $\uparrow$ & $100$ & $200$ & & $100$ & $1200$ \\
& $\downarrow$ & $50$ & $250$ & & $100$ & $1200$ \\
\multirow{2}{*}{$tc_6$} & $\uparrow$ & $100$ & $200$ & & $40$ & $200$ \\
& $\downarrow$ & $50$ & $250$ & & $80$ & $200$ \\
\multirow{2}{*}{$tc_7$} & $\uparrow$ & $100$ & $200$ & & $40$ & $400$ \\
& $\downarrow$ & $50$ & $250$ & & $-$ & $-$ \\
\multirow{2}{*}{$tc_8$} & $\uparrow$ & $100$ & $100$ & & $-$ & $-$ \\
& $\downarrow$ & $100$ & $100$ & & $80$ & $400$ \\
\multirow{2}{*}{$tc_9$} & $\uparrow$ & $100$ & $300$ & & $40$ & $200$ \\
& $\downarrow$ & $-$ & $-$ & & $80$ & $200$ \\
\multirow{2}{*}{$tc_{10}$} & $\uparrow$ & $-$ & $-$ & & $100$ & $1200$ \\
& $\downarrow$ & $100$ & $300$ & & $100$ & $1200$ \\
\multirow{2}{*}{$tc_{11}$} & $\uparrow$ & $50$ & $250$ & & $100$ & $1250$ \\
& $\downarrow$ & $100$ & $200$ & & $100$ & $1250$ \\
\multirow{2}{*}{$tc_{12}$} & $\uparrow$ & $200$ & $100$ & & $80$ & $1200$ \\
& $\downarrow$ & $50$ & $250$ & & $80$ & $1200$ \\
\multirow{2}{*}{$tc_{13}$} & $\uparrow$ & $50$ & $250$ & & $40$ & $200$ \\
& $\downarrow$ & $100$ & $200$ & & $80$ & $200$ \\
\multirow{2}{*}{$tc_{14}$} & $\uparrow$ & $100$ & $50$ & & $40$ & $400$ \\
& $\downarrow$ & $-$ & $-$ & & $-$ & $-$ \\
\multirow{2}{*}{$tc_{15}$} & $\uparrow$ & $-$ & $-$ & & $-$ & $-$ \\
& $\downarrow$ & $100$ & $50$ & & $80$ & $400$ \\
\multirow{2}{*}{$tc_{16}$} & $\uparrow$ & $100$ & $100$ & & $40$ & $200$ \\
& $\downarrow$ & $-$ & $-$ & & $80$ & $200$ \\
\multirow{2}{*}{$tc_{17}$} & $\uparrow$ & $-$ & $-$ & & $100$ & $100$ \\
& $\downarrow$ & $100$ & $100$ & & $100$ & $100$ \\
\bottomrule
\end{tabular}
\end{table}

\clearpage

\bibliographystyle{IEEEtran}
\bibliography{References}

\end{document}